%
%
%
\newif\ifsect\newif\iffinal
\secttrue\finalfalse
\def\thm #1: #2{\medbreak\noindent{\bf #1:}\if(#2\thmp\else\thmn#2\fi}
\def\thmp #1) { (#1)\thmn{}}
\def\thmn#1#2\par{\enspace{\sl #1#2}\par
        \ifdim\lastskip<\medskipamount \removelastskip\penalty 55\medskip\fi}
\def\square{{\msam\char"03}}
\def\qedn{\thinspace\null\nobreak\hfill\square\par\medbreak}
\def\pf{\ifdim\lastskip<\smallskipamount \removelastskip\smallskip\fi
        \noindent{\sl Proof\/}:\enspace}
\def\itm#1{\item{\rm #1}\ignorespaces}

\def\bar#1{\overline{#1}}
%
%
%
%
\newcount\parano
\newcount\eqnumbo
\newcount\thmno
\newcount\versiono
\newbox\notaautore
\def\neweqt#1$${\xdef #1{(\number\parano.\number\eqnumbo)}
    \eqno #1$$
    \global \advance \eqnumbo by 1}
\def\newrem #1\par{\global \advance \thmno by 1
    \medbreak
{\bf Remark \the\parano.\the\thmno:}\enspace #1\par
\ifdim\lastskip<\medskipamount \removelastskip\penalty 55\medskip\fi}
\def\newthmt#1 #2: #3{ \global \advance \thmno by 1\xdef #2{\number\parano.\number\thmno}
    \medbreak\noindent
    {\bf #1 #2:}\if(#3\thmp\else\thmn#3\fi}
\def\neweqf#1$${\xdef #1{(\number\eqnumbo)}
    \eqno #1$$
    \global \advance \eqnumbo by 1}
\def\newthmf#1 #2: #3{    \global \advance \thmno by 1\xdef #2{\number\thmno}
    \medbreak\noindent
    {\bf #1 #2:}\if(#3\thmp\else\thmn#3\fi}
\def\forclose#1{\hfil\llap{$#1$}\hfilneg}
\def\newforclose#1{
	\ifsect\xdef #1{(\number\parano.\number\eqnumbo)}\else
	\xdef #1{(\number\eqnumbo)}\fi
	\hfil\llap{$#1$}\hfilneg
	\global \advance \eqnumbo by 1
	\iffinal\else\rsimb#1\fi}
\def\forevery#1#2$${\displaylines{\let\eqno=\forclose
        \let\neweq=\newforclose\hfilneg\rlap{$\qquad\quad\forall#1$}\hfil#2\cr}$$}
\def\noNota #1\par{}
\def\today{\ifcase\month\or
   January\or February\or March\or April\or May\or June\or July\or August\or
   September\or October\or November\or December\fi
   \space\number\year}
\def\inizia{\ifsect\let\neweq=\neweqt\else\let\neweq=\neweqf\fi
\ifsect\let\newthm=\newthmt\else\let\newthm=\newthmf\fi}
\def\bititolo{\empty}
\gdef\begin #1 #2\par{\xdef\titolo{#2}
\ifsect\let\neweq=\neweqt\else\let\neweq=\neweqf\fi
\ifsect\let\newthm=\newthmt\else\let\newthm=\newthmf\fi
\iffinal\let\Nota=\noNota\fi
\centerline{\titlefont\titolo}
\if\bititolo\empty\else\medskip\centerline{\titlefont\bititolo}
\xdef\titolo{\titolo\ \bititolo}\fi
\bigskip
\centerline{\bigfont
\autore \ifvoid\notaautore\else\footnote{${}^1$}{\unhbox\notaautore}\fi}
\bigskip\if\istituto!\centerline{\today}\else
\centerline{\istituto}
\centerline{\indirizzo}
\centerline{\email}
\medskip
\centerline{#1\ \anno}\fi
\bigskip\bigskip
\ifsect\else\global\thmno=1\global\eqnumbo=1\fi}
\def\anno{2012}
\def\raggedleft{\leftskip2cm plus1fill \spaceskip.3333em \xspaceskip.5em
\parindent=0pt\relax}

\font\titlefont=cmssbx10 scaled \magstep1
\font\bigfont=cmr12
\font\eightrm=cmr8
\font\sc=cmcsc10
\font\bbr=msbm10
\font\sbbr=msbm7
\font\ssbbr=msbm5
\font\msam=msam10

\font\bfm=cmmib10

\font\SS=cmss10

\nopagenumbers
\binoppenalty=10000
\relpenalty=10000
\newfam\amsfam
\textfont\amsfam=\bbr \scriptfont\amsfam=\sbbr \scriptscriptfont\amsfam=\ssbbr
\newfam\boldifam
\textfont\boldifam=\bfm
\let\de=\partial

\def\Hol{\mathop{\rm Hol}\nolimits}

\def\Fix{\mathop{\rm Fix}\nolimits}

\def\Rea{\mathop{\rm Re}\nolimits}

\def\bigoperp{\mathop{\hbox{$\bigcirc\kern-11.8pt\perp$}}\limits}

\mathchardef\void="083F
\mathchardef\ellb="0960
\mathchardef\taub="091C
\def\C{{\mathchoice{\hbox{\bbr C}}{\hbox{\bbr C}}{\hbox{\sbbr C}}
{\hbox{\sbbr C}}}}
\def\R{{\mathchoice{\hbox{\bbr R}}{\hbox{\bbr R}}{\hbox{\sbbr R}}
{\hbox{\sbbr R}}}}
\def\N{{\mathchoice{\hbox{\bbr N}}{\hbox{\bbr N}}{\hbox{\sbbr N}}
{\hbox{\sbbr N}}}}

\newcount\notitle
\notitle=1
\headline={\ifodd\pageno\rhead\else\lhead\fi}
\def\rhead{\ifnum\pageno=\notitle\iffinal\hfill\else\hfill\tt Version
\the\versiono; \the\day/\the\month/\the\year\fi\else\hfill\eightrm\titolo\hfill
\folio\fi}
\def\lhead{\ifnum\pageno=\notitle\hfill\else\eightrm\folio\hfill\autore\hfill
\fi}
\newbox\bibliobox
\def\setref #1{\setbox\bibliobox=\hbox{[#1]\enspace}
    \parindent=\wd\bibliobox}
\def\biblap#1{\noindent\hang\rlap{[#1]\enspace}\indent\ignorespaces}
\def\art#1 #2: #3! #4! #5 #6 #7-#8 \par{\biblap{#1}#2: {\sl #3\/}.
    #4 {\bf #5} (#6)\if.#7\else, \hbox{#7--#8}\fi.\par\smallskip}
\def\book#1 #2: #3! #4 \par{\biblap{#1}#2: {\bf #3.} #4.\par\smallskip}
\def\coll#1 #2: #3! #4! #5 \par{\biblap{#1}#2: {\sl #3\/}. In {\bf #4,}
#5.\par\smallskip}
\def\pre#1 #2: #3! #4! #5 \par{\biblap{#1}#2: {\sl #3\/}. #4, #5.\par\smallskip}
\def\Bittersweet{\pdfsetcolor{0. 0.75 1. 0.24}} 
\def\Black{\pdfsetcolor{0. 0. 0. 1.}} 
\def\raggedleft{\leftskip2cm plus1fill \spaceskip.3333em \xspaceskip.5em
\parindent=0pt\relax}
\def\Nota #1\par{\medbreak\begingroup\Bittersweet\raggedleft
#1\par\endgroup\Black
\ifdim\lastskip<\medskipamount \removelastskip\penalty 55\medskip\fi}
%

%
\def\smallsect #1. #2\par{\bigbreak\noindent{\bf #1.}\enspace{\bf #2}\par
    \global\parano=#1\global\eqnumbo=1\global\thmno=0
    \nobreak\smallskip\nobreak\noindent\message{#2}}
\def\newdef #1\par{\global \advance \thmno by 1
    \medbreak
{\bf Definition \the\parano.\the\thmno:}\enspace #1\par
\ifdim\lastskip<\medskipamount \removelastskip\penalty 55\medskip\fi}
\finaltrue
\versiono=8
%



\let\te=\textstyle

\def\mlog{\mathop{{\te{1\over2}\log}}}

\def\autore{Marco Abate${}^1$\footnote{}{\eightrm 2010 Mathematics Subject Classification: Primary 32H50, 37F99.\hfill\break\indent Keywords: Wolff-Denjoy theorem; convex domains; holomorphic dynamics.}, Jasmin Raissy${}^2$\footnote{${}^*$}{\eightrm Partially supported by the PRIN2009 grant ``Critical Point Theory and Perturbative Methods for Nonlinear Differential Equations".}}
\def\indirizzo{\vbox{\hfill${}^2$Institut de Math\'ematiques de Toulouse,
Universit\'e Paul Sabatier,\hfill\break\null\hfill
118 route de Narbonne, 31062 Toulouse, France.
E-mail: jraissy@math.univ-tlse.fr\hfill\null}}
\def\istituto{\vbox{\hfill${}^1$Dipartimento di Matematica, Universit\`a
di Pisa,\hfill\break\null\hfill Largo Pontecorvo 5, 56127 Pisa,
Italy. E-mail: abate@dm.unipi.it\hfill\null\vskip5pt}}
\def\email{}
\begin {October} Wolff-Denjoy theorems in non-smooth convex domains

{{\narrower\noindent{\sc Abstract.}  
We give a short proof of Wolff-Denjoy theorem for (not necessarily smooth) strictly convex domains. With similar techniques we are also able to prove a Wolff-Denjoy theorem for weakly convex domains, again without any smoothness assumption on the boundary.

}}

\smallsect 0. Introduction

Studying the dynamics of a holomorphic self-map $f\colon\Delta\to\Delta$ of the unit disk~$\Delta\subset\C$ one is naturally led to consider two different cases. If $f$ has a fixed point then Schwarz's lemma readily implies that either $f$ is an elliptic automorphism, or the sequence $\{f^k\}$ of iterates of~$f$ converges (uniformly on compact sets) to the fixed point. 
\def\autore{Marco Abate and Jasmin Raissy}
The classical Wolff-Denjoy theorem ([W], [D]) says what happens when $f$ has no fixed points:

\newthm Theorem \WD: (Wolff-Denjoy) Let $f\colon\Delta\to\Delta$ be a holomorphic self-map without fixed points. Then there exists a point $\tau\in\de\Delta$ such that the sequence $\{f^k\}$ of iterates of~$f$ converges (uniformly on compact sets) to the constant map~$\tau$.

Since its discovery, a lot of work has been devoted to obtain similar statements in more general situations (surveys covering different aspects of this topic are [A3, RS, ES]). In one complex variable, there are results in multiply connected domains, multiply and infinitely connected Riemann surfaces, and even in the settings of one-parameter semigroups and of random dynamical systems (see, e.g., [H, L, B]). In several complex variables, the first Wolff-Denjoy theorems are due to Herv\'e [He1, 2];
in particular, in [He2] he proved a statement identical to the one above for fixed points free self-maps of the unit ball $B^n\subset\C^n$. Herv\'e's theorem has also been generalized in various ways to open unit balls of complex Hilbert and Banach spaces (see, e.g., [BKS, S] and references therein). 

A breakthrough occurred in 1988, when the first author (see [A1]) showed how to prove a Wolff-Denjoy theorem for holomorphic self-maps of smoothly bounded strongly convex domains in~$\C^n$. The techniques introduced there turned out to be quite effective in other contexts too (see, e.g., [A5, AR, Br1, Br2]); but in particular they led to Wolff-Denjoy theorems in smooth strongly pseudoconvex domains and smooth domains of finite type (see, e.g., [A4, Hu, RZ, Br3]). 

Two natural questions were left open by the previous results : how much does the boundary smoothness matter? And, what happens in weakly (pseudo)convex domains? As already shown by the results obtained by Herv\'e [He1] in the bidisk, if we drop both boundary smoothness and strong convexity the situation becomes much more complicated; but most of Herv\'e's techniques
were specific for the bidisk, and so not necessarily applicable to more general domains.
On the other hand, for smooth weakly convex domains a Wolff-Denjoy theorem was already obtained in [A3] (but here we shall get a better result; see Corollary~3.2). 

In 2012, Budzy\'nska [Bu2] (see also [BKR] and [Bu3] for infinite dimensional generalizations) finally proved a Wolff-Denjoy theorem for holomorphic fixed point free self-maps of a bounded strictly convex domain in~$\C^n$, under no smoothness assumption on the boundary; but she did not deal with weakly convex domains.

In Section 2 of this paper (Section 1 is devoted to recalling a few known preliminary facts) we shall give a simpler proof of Budzy\'nska's result, using only tools already introduced in [A1] and no additional machinery; it is worth mentioning that the final proof is also simpler than the proof presented in [A1] for the smooth case. In Section 3 we shall furthermore show how, combining our ideas with Budzy\'nska's new tools, one can obtain a Wolff-Denjoy theorem for weakly convex domains with no smoothness assumptions, thus addressing the second natural question mentioned above. Finally, in Section 4 we shall specialize our results to the polydisk, and we shall see that Herv\'e's results imply that our statements are essentially optimal. 
 
\smallsect 1. Preliminaries

In this section we shall collect a few more or less known facts on bounded convex domains in $\C^n$. 

\smallskip
\noindent{\it 1.1. Euclidean geometry.}

\smallskip\noindent
Let us begin by recalling a few standard 
definitions and notations.

\newdef Given $x$,~$y\in\C^n$ let
$$
[x,y]=\{sx+(1-s)y\in\C^n\mid s\in[0,1]\}\qquad\hbox{and}\qquad
(x,y)=\{sx+(1-s)y\in\C^n\mid s\in(0,1)\}
$$
denote the {\sl closed,} respectively {\sl open, segment} connecting $x$ and $y$. 
A set $D\subseteq\C^n$ is {\sl convex} if $[x,y]\subseteq D$ for all $x$,~$y\in D$;
and {\sl strictly convex} if $(x,y)\subseteq D$ for all $x$,~$y\in\bar{D}$.

An easy but useful observation is:

\newthm Lemma \couno: Let $D\subset\C^n$ be a convex domain. Then:
\itm{(i)} $(z,w)\subset D$ for all $z\in D$ and $w\in\de D$;
\itm{(ii)} if $x$, $y\in\de D$ then either $(x,y)\subset\de D$ or $(x,y)\subset D$.

%

This suggests the following

\newdef Let $D\subset\C^n$ be a convex domain. Given $x\in\de D$, we put
$$
{\rm ch}(x)=\{y\in\de D\mid [x,y]\subset\de D\}\;;
$$
we shall say that $x$ is a {\sl strictly convex point} if ${\rm ch}(x)=\{x\}$. More generally,
given $F\subseteq\de D$ we put
$$
{\rm ch}(F)=\bigcup_{x\in F}{\rm ch}(x)\;.
$$

A similar construction having a more holomorphic character is the following:

\newdef Let $D\subset\C^n$ be a convex domain. A {\sl complex supporting functional} at~$x\in\de D$ is a
$\C$-linear map $\sigma\colon\C^n\to\C$ such that $\Rea\sigma(z)<\Rea\sigma(x)$ for all $z\in D$. A {\sl complex supporting hyperplane} at $x\in\de D$ is 
an affine complex hyperplane $L\subset\C^n$ of the form
$L=x+\ker\sigma$, where $\sigma$ is a complex supporting functional at~$x$ (the existence of complex supporting functionals and hyperplanes is guaranteed by the Hahn-Banach theorem).
Given $x\in\de D$, we shall denote by
${\rm Ch}(x)$ the intersection of~$\bar{D}$ with of all complex supporting hyperplanes at~$x$. Clearly, ${\rm Ch}(x)$ is a closed convex set containing~$x$; in particular, ${\rm Ch}(x)\subseteq{\rm ch}(x)$. If ${\rm Ch}(x)=\{x\}$ we say that $x$ is a {\sl strictly $\C$-linearly convex point;}
and we say that $D$ is {\sl strictly $\C$-linearly convex} if all points of~$\de D$ are
strictly $\C$-linearly convex. Finally, if $F\subset\de D$ we set
$$
{\rm Ch}(F)=\bigcup_{x\in F}{\rm Ch}(x)\subseteq{\rm ch}(F)\;.
$$


\newrem If $\de D$ is of class $C^1$ then for each $x\in\de D$ there exists a unique complex
supporting hyperplane at~$x$, and thus ${\rm Ch}(x)$ coincides with the intersection of the complex supporting hyperplane with~$\de D$, which is smaller than the
flat region introduced in [A3, p. 277] as the intersection of $\de D$ with the real supporting hyperplane. But non-smooth points can have more than one complex supporting hyperplanes; this
happens for instance in the polydisk (see Section~4).  

\bigbreak

\noindent{\it 1.2. Intrinsic geometry}
\smallskip
\noindent
The intrinsic (complex) geometry of convex domains is conveniently described using the (intrinsic) Kobayashi distance. We refer to [A3], [JP] and [K] for details and much more on the Kobayashi (pseudo)distance in complex manifolds; here we shall just recall what is needed for our aims. Let $k_\Delta$ denote the 
Poincar\'e distance on the unit disk~$\Delta\subset{\bf C}$. If $X$ is a complex manifold,
the {\sl Lempert function}~$\delta_X\colon X\times X\to{\bf R}^+$ of~$X$ is 
$$
\delta_X(z,w)=\inf\{k_\Delta(\zeta,\eta)\mid
\hbox{there exists a holomorphic map $\phi\colon\Delta\to X$ with $\phi(\zeta)=z$ and
$\phi(\eta)=w$}\}
$$
for all $z$, $w\in X$. In general, the {\sl Kobayashi pseudodistance}~$k_X\colon X\times X\to{\bf R}^+$ of~$X$
is the largest pseudodistance on~$X$ bounded above by~$\delta_X$; when $D\subset\subset\C^n$ is a bounded convex domain in~$\C^n$, Lempert [Le] has proved that $\delta_D$ is an actual distance, and thus it coincides with the Kobayashi distance~$k_D$ of~$D$.

The main property of the Kobayashi (pseudo)distance is that it is contracted by holomorphic maps: if $f\colon X\to Y$ is a holomorphic map then
$k_Y\bigl(f(z),f(w)\bigr)\le k_X(z,w)$ for all $z$,~$w\in X$.
In particular, biholomorphisms are isometries, and holomorphic self-maps are $k_X$-nonexpansive.

The Kobayashi distance of convex domains enjoys several interesting properties; for instance,
it coincides with the Carath\'eodory distance, and it is 
a complete distance (see, e.g., [A3] or [Le]); in particular, $k_D$-bounded subsets of~$D$ are relatively
compact in~$D$. We shall also need
is the following estimates:

\newthm Lemma \codue: ([Le, KKR1, KS]) Let $D\subset\subset\C^n$ be a bounded convex domain. Then:
\itm{(i)} if $z_1$, $z_2$, $w_1$, $w_2\in D$ and $s\in [0, 1]$ then
$$
k_D\bigl(sz_1 + (1-s)w_1, sz_2 + (1-s)w_2\bigr)\le \max\bigl\{k_D(z_1,z_2), k_D(w_1, w_2)\bigr\}\;;
$$
\itm{(ii)} if $z$, $w\in D$ and $s$, $t\in[0,1]$ then
$$
k_D\bigl(sz + (1- s)w, tz + (1-t )w\bigr)\le k_D (z, w)\;.
$$

%
%

As a consequence we have:

\newthm Lemma \duno: Let $D\subset\subset\C^n$ be a bounded convex domain, $x$, $y\in\de D$, and let $\{z_\nu\}$, $\{w_\nu\}\subset D$ be two sequences converging to $x$ and $y$ respectively.
If
$$
\sup_{\nu\in\N} k_D (z_\nu, w_\nu) = c< +\infty
$$ 
then $[x,y]\subset\de D$. In particular, if $x$ (or $y$) is a strictly convex point then $x=y$.

\pf By Lemma~\couno\ we know that either $(x, y) \subset D$, or $(x, y) \subset \de D$. Assume by contradiction that $(x, y) \subset D$. 
Lemma~\codue\ yields
$$
k_D\bigl(s z_\nu + (1-s)w_\nu, t z_\nu + (1-t) w_\nu\bigr)\le k_D(z_\nu, w_\nu) \le c
$$
for each $\nu\in\N$ and for all $s, t\in(0,1)$. Hence
$$
k_D\bigl(s x + (1-s)y, tx + (1-t) y\bigr)= \lim_{\nu\to\infty} k_D\bigl(s z_\nu + (1-s)w_\nu, t z_\nu + (1-t) w_\nu\bigr) \le c
$$
for all $s$, $t\in(0,1)$. But this implies that $(x,y)$ is relatively compact in~$D$, which is impossible because $x$,~$y\in\de D$.
\qedn

\smallskip
\noindent{\it 1.3. Dynamics}
\smallskip
\noindent In this subsection we recall a few known facts about the dynamics of holomorphic
(or more generally $k_D$-nonexpansive) self-maps of convex domains.

When $D\subset\subset\C^n$ is a bounded domain, by Montel's theorem the space
$\Hol(D,D)$ is relatively compact in~$\Hol(D,\C^n)$. In particular, if $f\in\Hol(D,D)$ then 
every sequence $\{f^{k_j}\}$ of iterates contains a subsequence converging to a 
holomorphic map $h\in\Hol(D,\C^n)$. Analogously, using this time Ascoli-Arzel\`a theorem,
if $f\colon D\to D$ is $k_D$-nonexpansive then every sequence $\{f^{k_j}\}$ of iterates contains a subsequence converging to a continuous map $h\colon D\to\bar{D}\subset\C^n$.

\newdef Let $D\subset\subset\C^n$ be a bounded domain, and $f\colon D\to D$ a holomorphic 
or $k_D$-nonexpansive self-map. A map $h\colon D\to\C^n$ is a {\sl limit point} of the sequence
$\{f^k\}$ of iterates of~$f$ if there is a subsequence $\{f^{k_j}\}$ of iterates converging (uniformly on compact subsets) to~$h$; we shall denote by $\Gamma(f)$ the set of all limit points of~$\{f^k\}$. The {\sl target set}~$T(f)$ of~$f$ is then defined as the union of the images of limit points of the sequence of iterates:
$$
T(f)=\bigcup_{h\in\Gamma(f)} h(D)\;.
$$

\newdef A sequence $\{f_k\}\subset C(X,Y)$ of continuous maps between topological spaces
is {\sl compactly divergent} if for each pair of compact subsets $H\subseteq X$ and $K\subseteq Y$
there is $k_0\in\N$ such that $f^k(H)\cap K=\void$ for all $k\ge k_0$. 

When $D$ is a convex domain, the target set either is contained
in~$D$ if $f$ has a fixed point or is contained in~$\de D$ if $f$ has no fixed points. More precisely,
we have the following statement (see [A1, A4, C, KKR1, KS, Bu1]):

\newthm Theorem \noncdiv: Let $D\subset\subset\C^n$ be a bounded convex domain,
and $f\colon D\to D$ a $k_D$-nonexpansive (e.g., holomorphic) self-map. Then the following
assertions are equivalent:
\smallskip
\itm{(i)} $f$ has a fixed point in $D$;
\itm{(ii)} the sequence $\{f^k\}$ is not compactly divergent;
\itm{(iii)} the sequence $\{f^k\}$ has no compactly divergent subsequences;
\itm{(iv)} $\{f^k(z)\}$ is relatively compact in~$D$ for all $z\in D$;
\itm{(v)} there exists $z_0\in D$ such that $\{f^k(z_0)\}$ is relatively compact in~$D$;
\itm{(vi)} there exists $z_0\in D$ such that $\{f^k(z_0)\}$ admits a subsequence 
relatively compact in $D$.

\newrem For more general taut domains (and $f$ holomorphic) the statements (ii)--(vi)
are still equivalent. For some classes of domains, these statements
are equivalent to $f$ having a periodic point (see [A4] and [Hu]); however, there exist
holomorphic self-maps of a taut topologically 
contractible smooth domain satisfying (ii)--(vi) but without fixed points
(see [AH]).

When the sequence of iterates of $f$ is not compactly divergent (i.e., when $f$ 
has a fixed point if $D$ is convex) then the target set of $f$ has already been characterized ([Be], [A1, 5]). 
In particular, using Theorem~\noncdiv\ and repeating word by word the proof of
[A3, Theorem~2.1.29] we obtain

\newthm Theorem \tsncd: Let $D\subset\subset\C^n$ be a bounded convex domain,
and $f\colon D\to D$ a $k_D$-nonexpansive (e.g., holomorphic) self-map of~$D$. Assume that $f$ has a fixed point in~$D$.
Then $T(f)$ is a $k_D$-nonexpansive (respectively, holomorphic) retract of~$D$.
More precisely, there exists a unique $k_D$-nonexpansive (respectively, holomorphic)
retraction $\rho\colon D\to T(f)$ which is a limit point of~$\{f^k\}$, such that every limit point
of $\{f^k\}$ is of the form $\gamma\circ\rho$, where $\gamma\colon T(f)\to T(f)$ is a
(biholomorphic) invertible $k_D$-isometry, and $f|_{T(f)}$ is a (biholomorphic) invertible
$k_D$-isometry. 

In this paper we want to describe the target set of fixed points free self-maps of bounded convex domains.

\smallsect 2. Strictly convex domains 

Since Wolff's proof of the Wolff-Denjoy theorem [W], horospheres have been the main tool needed for
the study of the dynamics of fixed points free holomorphic self-maps. Let us recall the general definitions introduced in [A1, 3].

\newdef Let $D\subset\subset\C^n$ be a bounded domain, $z_0\in D$, $x\in\de D$ and
$R>0$. The {\sl small horosphere} $E_{z_0}(x,R)$ and the {\sl large horosphere}
$F_{z_0}(x,R)$ of {\sl center}~$x$, {\sl pole}~$z_0$ and {\sl radius}~$R$ are defined by
$$
\eqalign{
E_{z_0}(x,R)&=\bigl\{z\in D\bigm| \limsup_{w\to x}\bigl[k_D(z,w)-k_D(z_0,w)\bigr]<
{\textstyle{1\over2}}\log R\bigr\}\;,\cr
F_{z_0}(x,R)&=\bigl\{z\in D\bigm| \liminf_{w\to x}\bigl[k_D(z,w)-k_D(z_0,w)\bigr]<
{\textstyle{1\over2}}\log R\bigr\}\;.\cr}
$$

The following lemma contains some basic properties of horospheres, immediate consequence of the definition and of Lemma~\codue:

\newthm Lemma \dimm: Let $D\subset\subset\C^n$ be a bounded domain, $z_0\in D$ and $x\in\de D$. Then:
\itm{(i)} we have $E_{z_0}(x,R)\subseteq F_{z_0}(x,R)$ for every $R>0$;
\itm{(ii)} we have $\bar{E_{z_0}(x,R_1)}\cap D\subseteq E_{z_0}(x,R_2)$ and $\bar{F_{z_0}(x,R_1)}\cap D\subseteq F_{z_0}(x,R_2)$ for every $0<R_1<R_2$;
\itm{(iii)} we have
$B_D\bigl(z_0,{1\over2}\log R\bigr)\subseteq E_{z_0}(x,R)$ for all $R>1$, where 
$B_D(z_0,r)$ denotes the Kobayashi ball of center~$z_0$ and radius~$r$;
\itm{(iv)} we have $F_{z_0}(x,R)\cap B_D\bigl(z_0,-{1\over2}\log R)=\void$ for all $0<R<1$;
\itm{(v)} $\bigcup\limits_{R>0}E_{z_0}(x,R)=\bigcup\limits_{R>0}F_{z_0}(x,R)=D$ and
$\bigcap\limits_{R>0}E_{z_0}(x,R)=\bigcap\limits_{R>0}F_{z_0}(x,R)=\void$;
\itm{(vi)} if moreover $D$ is convex then $E_{z_0}(x,R)$ is convex for every $R>0$.


Large horospheres are not always convex, even if $D$ is convex; an example is given
by the horospheres in the polydisk (see Section~4). However, large horospheres in
convex domains are always star-shaped with respect to the center:

\newthm Lemma \uuno: Let $D\subset\subset\C^n$ be a bounded convex domain, $z_0\in D$, $R>0$
and $x\in\de D$. Then we have $[x,z]\subset \bar{F_{z_0}(x,R)}$ for 
all $z\in\bar{F_{z_0}(x,R)}$. In particular, $x$ always belongs to $\bar{F_{z_0}(x,R)}$.

\pf Given $z\in F_{z_0}(x,R)$, choose a sequence $\{x_\nu\}\subset D$ converging to~$x$
and such that the limit of $k_D(z,x_\nu)-k_D(z_0,x_\nu)$ exists and is less than~$\mlog R$. 
Given $0<s<1$, let $h_\nu^s\colon D\to D$ be 
$$
\forevery{w\in D}h_\nu^s(w)=sw+(1-s)x_\nu\;;
$$
then $h_\nu^s(x_\nu)=x_\nu$, and moreover
$$
\forevery{z_1,z_2\in D} k_D\bigl(h_\nu^s(z_1),h_\nu^s(z_2)\bigr)\le k_D(z_1,z_2)
$$
because $h_\nu^s$ is a holomorphic self-map of~$D$.
In particular,
$$
\limsup_{\nu\to+\infty}\bigl[k_D\bigl(h_\nu^s(z),x_\nu)-k_D(z_0,x_\nu)\bigr]
\le \lim_{\nu\to+\infty}\bigl[k_D(z,x_\nu)-k_D(z_0,x_\nu)\bigr]<\mlog R\;.
$$
Furthermore we have
$$
\bigl|k_D\bigl(sz+(1-s)x,x_\nu\bigr)-k_D\bigl(h^s_\nu(z),x_\nu\bigr)\bigr|
\le k_D\bigl(sz+(1-s)x_\nu,sz+(1-s)x\bigr)\to 0
$$
as $\nu\to+\infty$. Therefore
$$
\eqalign{
\liminf_{w\to x}&\bigl[k_D\bigl(sz+(1-s)x,w\bigr)-k_D(z_0,w)\bigr]
\le\limsup_{\nu\to+\infty}\bigl[k_D\bigl(sz+(1-s)x,x_\nu\bigr)-k_D(z_0,x_\nu)\bigr]\cr
&\le
\limsup_{\nu\to+\infty}\bigl[k_D\bigl(h_\nu^s(z),x_\nu\bigr)-k_D(z_0,x_\nu)\bigr]
+\lim_{\nu\to+\infty}\bigl[k_D\bigl(sz+(1-s)x,x_\nu\bigr)-k_D\bigl(h_\nu^s(z),x_\nu\bigr)\bigr]\cr
&<\mlog R\;,\cr}
$$
and thus $sz+(1-s)x\in F_{z_0}(x,R)$. Letting $s\to 1$ we get $x\in\bar{F_{z_0}(x,R)}$,
and we have proved the assertion for $z\in F_{z_0}(x,R)$. If $z\in\de F_{z_0}(x,R)$,
it suffices to apply the statement to a sequence in~$F_{z_0}(x,R)$ approaching~$z$.
\qedn

One of the main points in the proof given in [A1] of the Wolff-Denjoy theorem for strongly convex $C^2$ domains is the fact that in such domains the intersection between the closure of a large horosphere
and the boundary of the domain reduces to the center of the horosphere. The following
corollary will play the same r\^ole for not necessarily smooth convex domains:

\newthm Corollary \udue: Let $D\subset\subset\C^n$ be a bounded convex domain, $z_0\in D$, and $x\in\de D$. Then 
$$
\bigcap_{R>0} \bar{F_{z_0}(x,R)}\subseteq{\rm ch}(x)\;.
\neweq\equinter
$$
In particular, if $x$ is a strictly convex point then $\bigcap\limits_{R>0} \bar{F_{z_0}(x,R)}=\{x\}$.

\pf First of all, Lemma~\dimm\ implies that the intersection in \equinter\ is not empty and contained in~$\de D$.
Take $\tilde x\in\bigcap_{R>0} \bar{F_{z_0}(x,R)}$ different from~$x$. Then
Lemma~\uuno\ implies that the whole
segment $[x,\tilde x]$ is contained in the intersection, and thus in~$\de D$; hence $\tilde x\in{\rm ch}(x)$,
and we are done.
\qedn

Let us turn now to the study of the target set. A first step in this direction is the
following:

%
%

\newthm Proposition \AV: 
Let $D\subset\subset\C^n$ be a bounded convex domain. Then:
\smallskip
\itm{(i)} for every connected complex manifold $X$ and every holomorphic map
$h\colon X\to\C^n$ such that $h(X)\subset\bar{D}$ and $h(X)\cap\de D\ne\void$ we have 
$$
h(X)\subseteq\bigcap_{x\in X}{\rm Ch}\bigl(h(x)\bigr)\subseteq\de D\;.
$$
In particular, if $h$ is a limit point of the sequence of iterates of a holomorphic self-map of $D$
without fixed points we have
$$
h(D)\subseteq\bigcap_{z\in D}{\rm Ch}\bigl(h(z)\bigr)\;.
$$
\itm{(ii)}Let $f\colon D\to D$ be a $k_D$-nonexpansive self-map without fixed points, and $h\colon D\to\C^n$ a limit point of $\{f^k\}$. Then
$$
h(D)\subseteq \bigcap_{z\in D}\hbox{\rm ch}\bigl(h(z)\bigr)\;.
$$

\pf (i) The fact that $h(X)\subseteq\de D$ is an immediate consequence of the maximum
principle (see, e.g., [AV, Lemma~2.1]). 

Let now $L=h(x_0)+\ker\sigma$ be a complex supporting hyperplane at~$h(x_0)$. Then 
$\Rea(\sigma\circ h)\le\Rea\sigma\bigl(h(x_0)\bigr)$ on~$X$; therefore, by the maximum principle, $\sigma\circ h\equiv\sigma\bigl(h(x_0)\bigr)$, that is $h(X)\subset L$.
Since this holds for all complex supporting hyperplanes at~$h(x_0)$ the assertion follows.

(ii) Let $\{f^{k_j}\}$ be a subsequence of iterates converging to $h$. Since $f$ has no fixed points,
we know by Theorem~\noncdiv\ that $h(D)\subseteq\de D$. Furthermore
$$
\forevery{z,w\in D} k_D\bigl(f^{k_j}(z),f^{k_j}(w)\bigr)\le k_D(z,w)<+\infty\;;
$$
therefore Lemma~\duno\ implies $[h(z),h(w)]\subset\de D$, and the assertion follows.\qedn

The disadvantage of these statements is that the right-hand side still depends on the 
given limit point of the sequence of iterates; instead we would like to determine a 
subset of the boundary containing the whole target set. This can be accomplished as follows:

\newthm Lemma \dpuno: Let $D\subset\subset\C^n$ be a bounded convex domain, and $f\colon D\to D$ a $k_D$-nonexpansive (respectively, holomorphic) self-map without fixed points. Assume there exist $\void\ne E\subseteq F\subset D$ such that $f^k(E) \subset F$ for all~$k\in\N$. Then
we have
$$
T(f)\subseteq \hbox{\rm ch}\bigl(\bar F\cap \de D\bigr)
$$
if $f$ is $k_D$-nonexpansive, or 
$$
T(f)\subseteq \hbox{\rm Ch}\bigl(\bar F\cap \de D\bigr)
$$
if $f$ is holomorphic.

\pf Let $h$ be a limit point of the sequence of iterates of~$f$.
Since $f$ has no fixed points, we know that $h(D)\subseteq\de D$. Take $z_0\in E$; by 
assumption, the whole orbit of~$z_0$ is contained in~$F$. Therefore $h(z_0)\in\bar F\cap\de D$,
and the assertion follows from Proposition~\AV.
\qedn

The Wolff lemma [A1, Theorem~2.3], whose
proof can easily be adapted to the case of $k_D$-nonexpansive maps, is exactly what we need
to apply Lemma~\dpuno: 

\newthm Lemma \Wolffuno: Let $D\subset\subset\C^n$ be a convex domain, and let 
$f\colon D\to D$ be $k_D$-nonexpansive and without fixed points. Then there exists
$x\in\de D$ such that for every $z_0\in D$, $R>0$ and $k\in\N$ we have
$$
f^k\bigl(E_{z_0}(x,R)\bigr)\subseteq F_{z_0}(x,R)\;.
$$

%
%

We can now give a proof of a Wolff-Denjoy theorem for $k_D$-nonexpansive self-maps
of strictly convex domains in the same spirit as the proof given in [A1] of the
Wolff-Denjoy theorem for holomorphic self-maps of $C^2$ strongly convex domains,
without requiring the machinery introduced in [Bu2] and [BKR]:

\newthm Theorem \utre: Let $D\subset\subset\C^n$ be a bounded strictly convex domain,
and $f\colon D\to D$ a $k_D$-nonexpansive (e.g., holomorphic) self-map
without fixed points. Then there exists a $x_0\in\de D$ such that $T(f)=\{x_0\}$, that is
the sequence of iterates $\{f^k\}$ converges to the constant map~$x_0$.

\pf Fix $z_0\in D$. Lemmas~\Wolffuno\ and \dpuno\ give $x_0\in\de D$ such that
$$
T(f)\subseteq\bigcap_{R>0}{\rm ch}\bigl(\bar{F_{z_0}(x_0,R)}\cap\de D\bigr)\;.
$$
But $D$ is strictly convex; therefore ${\rm ch}\bigl(\bar{F_{z_0}(x_0,R)}\cap\de D\bigr)=
\bar{F_{z_0}(x_0,R)}\cap\de D$, and the assertion follows from Corollary~\udue.\qedn

In [A2] the first author characterized converging one-parameter semigroups of holomorphic self-maps
of smooth strongly convex domains. Theorem~\utre\ allows us to extend that characterization to (not necessarily smooth) strictly convex domains:

\newthm Corollary \inpiu: Let $D\subset\subset\C^n$ be a bounded strictly convex domain,
and $\Phi\colon\R^+\to\Hol(D,D)$ a one-parameter semigroup of holomorphic self-maps
of~$D$. Then $\Phi$ converges if and only if
\smallskip
\item{\rm(i)} either $\Phi$ has a fixed point $z_0\in D$ and the spectral generator at $z_0$ of $\Phi$
has no nonzero purely imaginary eigenvalues, or
\item{\rm(ii)} $\Phi$ has no fixed points.

\pf It follows arguing as in [A2, Theorem~1.3], replacing the references to [A1] by Theorem~\utre. \qedn

\smallsect 3. Weakly convex domains

As mentioned in the introduction, this approach works too
when $D$ is convex but not strictly convex. Simply applying the same argument used to prove Theorem~\utre\ one obtains
$$
T(f)\subseteq\bigcap_{R>0}{\rm ch}\bigl(\bar{F_{z_0}(x_0,R)}\cap\de D\bigr)
$$
in the $k_D$-nonexpansive case, and
$$
T(f)\subseteq\bigcap_{R>0}{\rm Ch}\bigl(\bar{F_{z_0}(x_0,R)}\cap\de D\bigr)
\neweq\eqhol
$$
in the holomorphic case (and it is easy to see that these intersections do not depend on~$z_0\in D$). This already can be used to strengthen the Wolff-Denjoy theorem obtained in [A3, Theorem~2.4.27] for weakly convex $C^2$ domains. Indeed, we can prove the following:

\newthm Proposition \aggiunta: Let $D\subset\subset\C^n$ be a $C^2$ bounded convex domain,
and $x\in\de D$. Then for every $z_0\in D$ and $R>0$ we have
$$
\bar{F_{z_0}(x,R)}\cap\de D\subseteq \hbox{\rm Ch}(x)\;.
$$
In particular, if $x$ is a strictly $\C$-linearly convex point then $\bar{F_{z_0}(x,R)}\cap\de D=\{x\}$.

\pf For every $x\in\de D$ let ${\bf n}_x$ denote the unit outer normal vector to~$\de D$ in~$x$,
and put $\sigma_x(z)=(z,{\bf n}_x)$,
where $(\cdot\,,\cdot)$ is the canonical Hermitian product. Then $\sigma_x$ is a complex supporting functional at~$x$ such that $\sigma_x(y)=\sigma_x(x)$ for some $y\in\de D$ if and only if $y\in\hbox{\rm Ch}(x)$. 

We can now argue as in the proof of [A3, Proposition~2.4.26] replacing the $P$-function
$\Psi\colon\de D\times\C^n\to\C$ given by $\Psi(x,z)=\exp\bigl(\sigma_x(z)-\sigma_x(x)\bigr)$,
with the $P$-function $\widehat\Psi\colon\de D\times\C^n\to\C$ given by 
$$
\widehat\Psi(x,z)={1\over 1-\bigl(\sigma_x(z)-\sigma_x(x)\bigr)}\;.
$$
\qedn

\newthm Corollary \aggiuntobis: Let $D\subset\subset\C^n$ be a $C^2$ bounded convex domain,
and $f\colon D\to D$ a holomorphic self-map without fixed points. Then there exists 
$x_0\in\de D$ such that $T(f)\subseteq\hbox{\rm Ch}(x_0)$. In particular, if $D$ is strictly
$\C$-linearly convex then the sequence of iterates $\{f^k\}$ converges to the constant map~$x_0$.

\pf It follows from \eqhol, Proposition~\aggiunta, and the fact that in $C^2$ convex domains each point in the boundary admits a unique complex supporting hyperplane.\qedn

\newrem We conjecture that the final assertion of this corollary should also hold for not necessarily smooth strictly $\C$-linearly convex domains. 

In weakly convex non-smooth domains large horospheres might be too large, and the right-hand side of \eqhol\ might coincide with the whole boundary of the domain (see Section~4 for an example in the polydisk); so to get an effective statement we need to replace them with smaller sets. 

Small horospheres might be too small; as shown by Frosini [F], there are holomorphic self-maps of the polydisk with no invariant
small horospheres. We thus need another kind of horospheres, defined by Kapeluszny, Kuczumow and Reich [KKR2], and studied in detail by Budzy\'nska [Bu2].
To introduce them we begin with a definition:

\newdef Let $D\subset\subset\C^n$ be a bounded domain, and $z_0\in D$. A sequence ${\bf x}=\{x_\nu\}\subset D$ converging to~$x\in\de D$ is
a {\sl horosphere sequence} at~$x$ if
the limit of $k_D(z,x_\nu)-k_D(z_0,x_\nu)$ as $\nu\to+\infty$ exists for all~$z\in D$.

\newrem It is easy to see that the notion of horosphere sequence does not depend on the point~$z_0$.

\newrem In [BKR] it is shown that every sequence in~$D$ converging to~$x\in\de D$
contains a subsequence which is a horosphere sequence at~$x$.
In strongly convex $C^3$ domains all sequences converging to a boundary point are horosphere sequences (see [A3, Theorem 2.6.47] and [BPT]); in Section~4 we shall give an explicit example of horosphere sequence in the polydisk.

\newdef Let $D\subset\subset\C^n$ be a bounded convex domain. Given $z_0\in D$,
let ${\bf x}$ be a horosphere sequence at~$x\in\de D$, and take $R>0$.
Then the {\sl sequence horosphere} $G_{z_0}(x,R,{\bf x})$ is defined as
$$
G_{z_0}(x,R,{\bf x})=\bigl\{z\in D\bigm|\lim_{\nu\to+\infty}\bigl[k_D(z,x_\nu)-k_D(z_0,x_\nu)\bigr]
<{\textstyle{1\over2}}\log R\bigr\}\;.
$$

The basic properties of sequence horospheres are contained in the following:

\newthm Proposition \newhor: ([KKR2, Bu2, BKR]) Let $D\subset\subset\C^n$ be a bounded convex domain. Fix $z_0\in D$, and let ${\bf x}=\{x_\nu\}\subset D$ be a horosphere sequence
at~$x\in\de D$ for $z_0$. Then:
\itm{(i)} $E_{z_0}(x,R)\subseteq G_{z_0}(x,R,{\bf x})\subseteq F_{z_0}(x,R)$ for all $R > 0$; 
\itm{(ii)} $G _{z_0}(x,R,{\bf x})$ is nonempty and convex for all $R>0$;
\itm{(iii)}$\bar{G_{z_0}(x,R_1,{\bf x})}\cap D\subset G_{z_0} (x,R_2, {\bf x})$ for all $0<R_1<R_2$;
\itm{(iv)} $B_D(z_0,{1\over2}\log R)\subset G_{z_0}(x, R, {\bf x})$ for all $R>1$;
\itm{(v)} $B_D(z_0,-{1\over2}\log R)\cap G_{z_0}(x, R, {\bf x})=\void$ for all $0<R<1$;
\itm{(vi)} $\bigcup\limits_{R>0}G_{z_0}(x,R,{\bf x})=D$ and
$\bigcap\limits_{R>0}G_{z_0}(x,R,{\bf x})=\void$.

\newrem If $\bf x$ is a horosphere sequence at~$x\in\de D$ then it is not difficult to check that the family $\{G_z(x,1,{\bf x})\}_{z\in D}$ and the family $\{G_{z_0}(x,R,{\bf x})\}_{R>0}$
with, $z_0\in D$ given, coincide.

It turns out that 
we can always find invariant sequence horospheres:

\newthm Lemma \Wolffdue: Let $D\subset\subset\C^n$ be a convex domain, and let 
$f\colon D\to D$ be $k_D$-nonexpansive and without fixed points. Then there exists
$x\in\de D$ and a horosphere sequence $\bf x$ at~$x$ such that
$$
f\bigl(G_{z_0}(x,R,{\bf x})\bigr)\subseteq G_{z_0}(x,R,{\bf x})
$$
for every $z_0\in D$ and $R>0$.

\pf Arguing as in the proof of [A1, Theorem~2.3] we can find a sequence $\{f_\nu\}$
of $k_D$-contractions with a unique fixed point $x_\nu\in D$ such that $f_\nu\to f$
and $x_\nu\to x\in\de D$ as $\nu\to+\infty$. 
Up to a subsequence, we can also assume (Remark~3.6) that ${\bf x}=\{x_\nu\}$ is a horosphere
sequence at~$x$.

Now, for every $z\in D$ we have
$$
\bigl|k_D\bigl(f(z),x_\nu\bigr)-k_D\bigl(f_\nu(z),x_\nu\bigr)\bigr|\le
k_D\bigl(f_\nu(z),f(z)\bigr)\to 0
$$ 
as $\nu\to+\infty$. Therefore if $z\in G_{z_0}(x,R,{\bf x})$ we get
$$
\eqalign{
\lim_{\nu\to+\infty}&\bigl[k_D\bigl(f(z),x_\nu\bigr)-k_D(z_0,x_\nu)\bigr]\cr
&\le\limsup_{\nu\to+\infty}\bigl[k_D\bigl(f_\nu(z),x_\nu\bigr)-k_D(z_0,x_\nu)\bigr]
+\limsup_{\nu\to+\infty}\bigl[k_D\bigl(f(z),x_\nu\bigr)-k_D\bigl(f_\nu(z),x_\nu\bigr)\bigr]\cr
&\le\lim_{\nu\to+\infty}\bigl[k_D(z,x_\nu)-k_D(z_0,x_\nu)\bigr]<\mlog R\cr}
$$
because $f_\nu(x_\nu)=x_\nu$ for all $\nu\in\N$, and we are done.\qedn


Putting everything together we can prove the following Wolff-Denjoy theorem
for (not necessarily strictly or smooth) convex domains:

\newthm Theorem \dpqua: Let $D\subset\subset\C^n$ be a bounded convex domain,
and $f\colon D\to D$ a $k_D$-nonexpansive (respectively, holomorphic) self-map without fixed points. Then there exist $x\in \de D$ and a horosphere
sequence ${\bf x}$ at~$x$ such that for any $z_0\in D$ we have
$$
T(f)\subseteq \bigcap_{z\in D}\hbox{\rm ch}\bigl(\bar{G_z(x,1,{\bf x})}\cap\de D\bigr)=\bigcap_{R>0}\hbox{\rm ch}\bigl(\bar{G_{z_0}(x,R,{\bf x})}\cap\de D\bigr)
$$
if $f$ is $k_D$-nonexpansive, or 
$$
T(f)\subseteq \bigcap_{z\in D}\hbox{\rm Ch}\bigl(\bar{G_z(x,1,{\bf x})}\cap\de D\bigr)=\bigcap_{R>0}\hbox{\rm Ch}\bigl(\bar{G_{z_0}(x,R,{\bf x})}\cap\de D\bigr)
$$
if $f$ is holomorphic.

\pf The equality of the intersections (and thus the independence of~$z_0$) is an immediate
consequence of Remark~3.9. Then the assertion follows from Lemmas~\Wolffdue\ and~\dpuno.\qedn

\newthm Corollary \dpslc: Let $D\subset\subset\C^n$ be a bounded strictly $\C$-linearly convex domain, and $f\colon D\to D$ a holomorphic self-map of $D$ without fixed points. Then there
exist $x\in \de D$ and a horosphere
sequence ${\bf x}$ at~$x$ such that for any $z_0\in D$ we have
$$
T(f)\subseteq \bigcap_{z\in D}\bar{G_z(x,1,{\bf x})}=\bigcap_{R>0}\bar{G_{z_0}(x,R,{\bf x})}\;.
$$

\pf It follows immediately from Theorem~\dpqua\ and the definition of strictly $\C$-linearly convex domain.\qedn

\smallsect 4. The polydisk


The polydisk $\Delta^n\subset\C^n$ is the unit ball for the norm
$$
\|z\| = \max\{|z_j| \mid j=1,\dots, n\}\;,
$$
and therefore
$$
\forevery{z,w\in \Delta^n} k_{\Delta^n}(z,w) = {1\over 2}\log{1+\|\gamma_z(w)\|\over 1-\|\gamma_z(w)\|}\;,
$$
where
$$
\gamma_z(w) = \left({w_1-z_1\over 1-\bar z_1 w_1},\dots, {w_n-z_n \over 1-\bar z_n w_n}\right)
$$
is an automorphism of the polydisk with $\gamma_z(z) = 0$.

Thanks to the homogeneity of $\Delta^n$, we can restrict ourselves to consider only horospheres with pole $z_0$ at the origin, and we have (see [A3, chapter 2.4.2] for detailed computations) the following description for horospheres with center $\xi\in\de \Delta^n$ and radius $R>0$:
$$
E_O(\xi, R) = \left\{z\in\Delta^n\biggm| \max_j\left\{\left.{|\xi_j-z_j|^2\over 1- |z_j|^2}\right| |\xi_j|=1\right \}< R\right\}
$$
and
$$
F_O(\xi, R) = \left\{z\in\Delta^n\biggm| \min_j\left\{\left.{|\xi_j-z_j|^2\over 1- |z_j|^2}\right| |\xi_j|=1\right \} < R\right\}.
$$
Moreover, $E_O(\xi, R)=F_O(\xi, R)$ if and only if $\xi$ has only one component of modulus $1$. 

%
Given $\xi\in\de\Delta^n$, a not difficult computation shows that
$$
{\rm ch}(\xi)=\bigcup_{|\xi_j|=1}\{\eta\in\de\Delta^n\mid \eta_j=\xi_j\}\qquad\hbox{and}\qquad
%
{\rm Ch}(\xi)=\bigcap_{|\xi_j|=1}\{\eta\in\de\Delta^n\mid \eta_j=\xi_j\}\;.
$$
This implies that in the polydisk large horospheres are too large to give a sensible Wolff-Denjoy theorem. Indeed we have ${\rm ch}\bigl(\bar{F_O(\xi,R)}\cap\de\Delta^n\bigr)={\rm Ch}\bigl(\bar{F_O(\xi,R)}\cap\de\Delta^n\bigr)=
\de\Delta^n$ if $\xi$ has at least two components of modulus~1, and
$$
{\rm ch}\bigl(\bar{F_O(\xi,R)}\cap\de\Delta^n\bigr)={\rm Ch}\bigl(\bar{F_O(\xi,R)}\cap\de\Delta^n\bigr)=
\de\Delta^n\setminus\{\eta\in\de\Delta^n\mid \eta_{j_0}\ne\xi_{j_0},\ |\eta_j|<1 \hbox{ for $j\ne j_0$}\}
$$
if $|\xi_{j_0}|=1$ and $|\xi_j|<1$ for $j\ne j_0$. 

Let us then compute the sequence horospheres. Fix a horosphere sequence ${\bf x}=\{x_\nu\}$ converging to~$\xi\in\de\Delta^n$. Arguing as in [A3, chapter 2.4.2], we arrive to the following 
$$
G_O(\xi,R,{\bf x})
=\left\{z\in\Delta^n\biggm| \max_j\left\{\left.{|\xi_j-z_j|^2\over 1- |z_j|^2}\lim_{\nu\to+ \infty}\min_h\left\{{1-|x_{\nu,h}|^2\over 1-|x_{\nu,j}|^2}\right\} \right| |\xi_j|=1\right \} < R\right\}\;.
$$
Since if $|\xi_j|=1$ we clearly have 
$$
\alpha_j:=\lim_{\nu\to+\infty}\min_h\left\{{1-|x_{\nu,h}|^2\over 1-|x_{\nu,j}|^2}\right\}\le 1\;,
$$
we get
$$
G_O(\xi,R,{\bf x})=\left\{z\in\Delta^n\biggm| \max_j\left\{\alpha_j{|\xi_j-z_j|^2\over 1- |z_j|^2}\biggm| |\xi_j|=1\right\}  < R\right\}.
$$
In other words, we can write $G_O(\xi,R,{\bf x})$ as a product
$$
G_O(\xi,R,{\bf x})=\prod_{j=1}^n E_j\;,
$$
where, denoting by $E^\Delta(\sigma,R)\subset\Delta$ the standard horocycle of center $\sigma\in\de\Delta$, pole the origin and radius~$R>0$, we put
$$
E_j=\cases{\Delta&if $|\xi_j|<1$,\cr
E^\Delta(\xi_j,R/\alpha_j)&if $|\xi_j|=1$.\cr}
$$
As a consequence, 
$$
{\rm ch}\bigl(\bar{G_O(\xi,R,{\bf x})}\cap\de\Delta^n\bigr)=
{\rm Ch}\bigl(\bar{G_O(\xi,R,{\bf x})}\cap\de\Delta^n\bigr)=
\bigcup_{j=1}^n \bar\Delta\times\cdots\times C_j(\xi)\times\cdots\times\bar\Delta\;,
$$
where
$$
C_j(\xi)=\cases{\{\xi_j\}&if $|\xi_j|=1$,\cr
\de\Delta& if $|\xi_j|<1$.\cr}
$$
Notice that the right-hand sides do not depend either on $R$ or on the horosphere sequence ${\bf x}$, but only on~$\xi$.


So Theorem~\dpqua\ in the polydisk assumes the following form:

\newthm Corollary \dpquapd: Let $f\colon \Delta^n\to \Delta^n$ be a $k_{\Delta^n}$-nonexpansive (e.g., holomorphic) self-map without fixed points. Then there exists $\xi\in \de\Delta^n$ such that
$$
T(f)\subseteq \bigcup_{j=1}^n \bar\Delta\times\cdots\times C_j(\xi)\times\cdots\times\bar\Delta\;,
\neweq\equnion
$$
where
$$
C_j(\xi)=\cases{\{\xi_j\}&if $|\xi_j|=1$,\cr
\de\Delta& if $|\xi_j|<1$.\cr}
$$

This is the best one can do, in the sense that while it is true that the image of a limit point of the sequence of iterates of~$f$ is always contained in just one of the set appearing in the right-hand side of~\equnion, it is impossible to determine a priori in which one on the basis of the point~$\xi$ only;
it is necessary to know something more about the map~$f$. Indeed, Herv\'e has proved the following:

\newthm Theorem \Herve: ([He1]) Let $F=(f,g)\colon\Delta^2\to\Delta^2$ be a holomorphic self-map of the bidisk, and write $f_w=f(\cdot,w)$ and $g_z=g(z,\cdot)$. Assume that $\Fix(F)=\void$. Then one and only one of the following cases occurs:
\smallskip
\item{\rm (0)} if $g(z,w)\equiv w$ (respectively, $f(z,w)\equiv z$) then the sequence of iterates of~$F$ converges uniformly on compact sets to $h(z,w)=(\sigma,w)$, where $\sigma$ is the common Wolff point of the $f_w$'s (respectively, to $h(z,w)=(z,\tau)$, where $\tau$ is the common Wolff point
of the $g_z$'s);
\item{\rm (1)} if $\Fix(f_w)=\void$ for all $w\in\Delta$ and $\Fix(g_z)=\{y(z)\}\subset\Delta$ for all $z\in\Delta$ (respectively, if $\Fix(f_w)=\{x(w)\}$ and $\Fix(g_z)=\void$) then 
$T(f)\subseteq\{\sigma\}\times\bar\Delta$, where $\sigma\in\de\Delta$ is the common Wolff point of the~$f_w$'s (respectively, $T(f)\subseteq\bar\Delta\times\{\tau\}$, where $\tau$ is the common Wolff point
of the $g_z$'s);
\item{\rm (2)} if $\Fix(f_w)=\void$ for all $w\in\Delta$ and $\Fix(g_z)=\void$ for all $z\in\Delta$ then either $T(f)\subseteq\{\sigma\}\times\bar\Delta$ or
$T(f)\subseteq\bar\Delta\times\{\tau\}$, where $\sigma\in\de\Delta$ is the common Wolff point of the~$f_w$'s, and $\tau\in\de\Delta$ is the common Wolff point of the~$g_z$;
\item{\rm (3)} if $\Fix(f_w)=\{x(w)\}\subset\Delta$ for all $w\in\Delta$ and $\Fix(g_z)=\{y(z)\}\subset\Delta$
for all $z\in\Delta$ then there are\break\indent $\sigma$,~$\tau\in\de D$ such that the sequence of iterates converges to the constant map $(\sigma,\tau)$.

We end this paper providing, as promised, an example of horosphere sequence.
Given $\xi\in\de\Delta^n$, put $x_\nu=(1-1/\nu)^{1/2}\xi$; we claim that ${\bf x}=\{x_\nu\}$ is a horosphere sequence. Indeed, arguing as in [A3, chapter~2.4.2] we see it suffices to show that
$$
\max_{|\xi_j|=1}\left\{\min_h\left\{{1-|x_{\nu,h}|^2\over 1-|x_{\nu,j}|^2}\right\}
{|1-\bar{z_j}x_{\nu,j}|^2\over 1-|z_j|^2}\right\}
$$
converges as $\nu\to+\infty$. But indeed with this choice of $x_\nu$ we have
$$
\eqalign{
\max_{|\xi_j|=1}\left\{\min_h\left\{{1-|x_{\nu,h}|^2\over 1-|x_{\nu,j}|^2}\right\}
{|1-\bar{z_j}x_{\nu,j}|^2\over 1-|z_j|^2}\right\}&=\max_{|\xi_j|=1}\left\{
\min_h\bigl\{\nu(1-|\xi_h|^2)+|\xi_h|^2\bigr\}{|1-\bar{z_j}x_{\nu,j}|^2\over 1-|z_j|^2}\right\}\cr
&=\max_{|\xi_j|=1}\left\{{|1-\bar{z_j}x_{\nu,j}|^2\over 1-|z_j|^2}\right\}\to
\max_{|\xi_j|=1}\left\{{|\xi_j-z_j|^2\over 1-|z_j|^2}\right\}\;.\cr}
$$
In particular, using this horosphere sequence one obtains $G_O(\xi,R,{\bf x})=E_O(\xi,R)$.

\setref{KKR2}
\beginsection References

\art A1 M. Abate: Horospheres and iterates of holomorphic maps! Math. Z.! 198 1988 225-238

\art A2 M. Abate: Converging semigroups of holomorphic maps! Rend. Acc. Naz. Lincei! 82
1988 223-227 

%
\book A3 M. Abate: Iteration theory of holomorphic maps on taut manifolds! Mediterranean 
Press, Cosenza, 1989. {\SS http://www.dm.unipi.it/$\sim$abate/libri/libriric/libriric.html}

\art A4 M. Abate: Iteration theory, compactly divergent sequences and commuting 
holomorphic maps! Ann. Scuola Norm. Sup. Pisa! 18 1991 167-191

%
\coll A5 M. Abate: Angular derivatives in several complex variables! Real methods in
complex and CR geometry! Eds. D. Zaitsev, G. Zampieri, Lect. Notes in Math.
1848, Springer, Berlin, 2004, pp. 1--47

\art AH M. Abate, P. Heinzner: Holomorphic actions on contractible domains 
without fixed points! Math. Z.! 211 1992 547-555

\art AR M. Abate, J. Raissy: Backward iteration in strongly convex domains! Adv. Math.! 228 2011 2837-2854

\art AV M. Abate, J.-P. Vigu\'e: Common fixed points in hyperbolic Riemann
surfaces and convex domains! Proc. Am. Math. Soc.! 112 1991 503-512

\art B  A.F. Beardon: Repeated compositions of analytic maps! 
Comput. Methods Funct. Theory! 1 2001 235-248 


\art Be E. Bedford: On the automorphism group of a Stein manifold! Math. Ann.! 266
1983 215-227

%

\art Br1 F. Bracci: Fixed points of commuting holomorphic mappings other than the Wolff point!
Trans. Amer. Math. Soc.! 355 2003 2569-2584 

\art Br2 F. Bracci: Dilatation and order of contact for holomorphic self-maps of strongly convex domains! Proc. London Math. Soc.! 86 2003 131-152

\art Br3 F. Bracci: A note on random holomorphic iteration in convex domains! Proc. Edinb. Math. Soc.! 51 2008 297-304

\art BPT F. Bracci, G. Patrizio, S. Trapani: The pluricomplex Poisson kernel for strongly convex domains! Trans. Amer. Math. Soc.! 361 2009 979-1005

\art Bu1 M. Budzy\'nska: Local uniform linear convexity with respect to the Kobayashi distance! Abstr. Appl. Anal.! 2003 2003 367-373 

\art Bu2 M. Budzy\'nska: The Denjoy-Wolff theorem in $\C^n$! Nonlinear Anal.! 75 2012 22-29 

\pre Bu3 M. Budzy\'nska: The Denjoy-Wolff theorem for condensing mappings in a bounded and strictly convex domain in a complex Banach space! Preprint! 2012 

\pre BKR M. Budzy\'nska, T. Kuczumow, S. Reich: Theorems of Denjoy-Wolff type! To appear in Ann. Mat. Pura Appl.!  2011 

%

\art BKS M. Budzy\'nska, T. Kuczumow, T. S\l odkowski: Total sets and semicontinuity of the Kobayashi distance! 
Nonlinear Anal.! 47 2001 2793-2803 

\art C A. Ca\l ka: On conditions under which isometries have bounded orbits! Colloq. Math.! 48 1984 219-227


\art D A. Denjoy: Sur l'it\'eration des fonctions analytiques! C.R. Acad. Sci. Paris!
182 1926 255-257

%
\book ES M. Elin, D. Shoikhet: Linearization models for complex dynamical systems. Topics in univalent functions, functional equations and semigroup theory! Birkh\"auser Verlag, Basel, 2010 
%

\coll F C. Frosini: Dynamics on bounded domains! The $p$-harmonic equation and recent advances in analysis! Contemp. Math., 370, Amer. Math. Soc., Providence, RI, 2005, pp. 99--117

%

\coll H M.H. Heins: A theorem of Wolff-Denjoy type!
 Complex analysis! 
 Birkh\"auser, Basel,  1988, pp. 81--86 
 
\art He1 M. Herv\'e: It\'eration des transformations analytiques dans le bicercle-unit\'e! Ann. Sci. \'Ec. Norm. Sup.! 71 1954 1-28 

\art He2 M. Herv\'e: Quelques propri\'et\'es des applications analytiques d'une
boule \`a $m$ dimensions dans elle-m\^eme! J. Math. Pures Appl.! 42 
1963 117-147 

\art Hu X.J. Huang: A non-degeneracy property of extremal mappings and iterates of
 holomorphic self-mappings!
 Ann. Scuola Norm. Sup. Pisa!  21 1994  399-419 

%

\book JP M. Jarnicki, P. Pflug: Invariant distances and metrics in complex analysis! Walter de Gruyter, Amsterdam, 1993

\art KKR1 J. Kapeluszny, T. Kuczumow, S. Reich: The Denjoy-Wolff theorem for condensing holomorphic mappings! J. Funct. Anal.! 167 1999 79-93 

\art KKR2 J. Kapeluszny, T. Kuczumow, S. Reich: The Denjoy-Wollf theorem in the open unit ball of a strictly convex Banach space! Adv. Math! 143 1999 111-123 



\book K S. Kobayashi: Hyperbolic complex spaces! Springer, Berlin, 1998

%
%


\art KS T. Kuczumow, A. Stachura: Iterates of holomorphic and $k_D$-nonexpansive mappings in convex domains in $\C^n$! Adv. Math.! 81 1990 90-98 

%
%
\art L F. L\'arusson: A Wolff-Denjoy theorem for infinitely connected Riemann surfaces!
 Proc. Amer. Math. Soc.!  124 1996  2745-2750 
 
\art Le L. Lempert: La m\'etrique de Kobayashi et la r\'epr\'esentation des domaines
sur la boule! Bull. Soc. Math. Fr.! 109 1981 427-474 

%
%
%
%
%

\book RS S. Reich, D. Shoikhet: Nonlinear semigroups, fixed points, and geometry of domains in Banach spaces!
Imperial College Press, London,  2005  

%
\art RZ F. Ren, W. Zhang: Dynamics on weakly pseudoconvex domains!
 Chinese Ann. Math. Ser. B!  16 1995 467-476 


\art S A. Stachura: Iterates of holomorphic self-maps of the unit ball in Hilbert space!
Proc. Amer. Math. Soc..! 93 1985 88-90 

%

\art W J. Wolff: Sur une g\'en\'eralisation d'un th\'eor\`eme de Schwarz! C.R. Acad. 
Sci. Paris! 182 1926 918-920


\bye